\documentclass[11pt,english]{revtex4}

\usepackage[T1]{fontenc}
\usepackage{babel}
\usepackage{graphicx,epsfig,amssymb,amsbsy}
\usepackage{color}
\usepackage{subeqn}
\usepackage{ulem}
\usepackage{todonotes,hyperref}
\newcommand{\bb}[1]{\left({#1}\right)}					
\newcommand{\sq}[1]{\left[#1\right]}						
\newcommand{\cc}[1]{\left\{#1\right\}}					
\newcommand{\ord}[1]{{\sf O}\bb{#1}}					

\newcommand{\eqblock}[2]{\begin{equation}\label{#1}\begin{array}{rcl}#2\end{array}\end{equation}}
\newcommand{\fref}[1]{figure~\ref{#1}}
\newcommand{\eref}[1]{equation (\ref{#1})}
\newcommand{\esref}[1]{system (\ref{#1})}
\newcommand{\sref}[1]{section~\ref{#1}}
\newcommand{\Fref}[1]{Figure~\ref{#1}}

\newcommand{\Erefs}[1]{Equations (\ref{#1})}

\newcommand{\eps}{\varepsilon}

\begin{document}

\title{Nonsmooth analogues of slow-fast dynamics -- pinching at a folded node}

\author{Mathieu Desroches}\address{INRIA Paris-Rocquencourt, SISYPHE Project Team, Domaine de Voluceau - BP 105, 78153 Le Chesnay cedex, France, Mathieu.Desroches@inria.fr}
\author{Mike R. Jeffrey}\address{University of Bristol, Department of Engineering Mathematics, Queen's Building, Bristol BS8 1TR, UK, Mike.Jeffrey@bristol.ac.uk}

\begin{abstract}
The folded node is a singularity associated with loss of normal hyperbolicity in systems where mixtures of slow and fast timescales arise due to singular perturbations. Canards are special solutions that reveal a counteractive feature of the local dynamics, namely trajectories that flow {\it from} attractive regions of space {\it into} repulsive regions. An alternative way to model switches between timescales is using piecewise-smooth differential equations. There is presently no adequate theory or method for relating slow-fast and piecewise-smooth models. Here we derive the analogous piecewise-smooth system for the folded node by {\it pinching} phase space to sharpen the switch between timescales. The corresponding piecewise-smooth system contains a so-called two-fold singularity, and exhibits the same topology and number of canards as the slow-fast system. Thus pinching provides a piecewise-smooth approximation to a slow-fast system. 
\end{abstract}

\maketitle


The purpose of this paper is to study the relation between discontinuities and singular perturbations in dynamical systems, by focusing on their singularities. We do this by forming a piecewise-smooth model of an important singular perturbation problem using the method of pinching introduced in \cite{jd10}. 
Pinching replaces a slow-fast timescale separation with a discontinuous switch, replacing Fenichel's slow manifolds \cite{f79} with Filippov's sliding surfaces \cite{f88}. In particular, it provides a piecewise-smooth model that captures, both qualitatively and quantitatively, the intricate dynamics of the smooth singular perturbation problem. The study of such singularities in discontinuous systems predates those in singularly perturbed systems \cite{f64,russian1,russian2}, while the latter has seen more progress in the study of oscillatory dynamics, see e.g. \cite{des12}. 
Here we show that the analogy between singularities in the two types of system is more than superficial, paving the way for a more rigorous study of piecewise-smooth systems as a means to studying singular perturbations in the future. 

\section{Introduction}\label{sec:intro}

Highly nonlinear changes in the dynamics of a system can be modelled by differential equations whose solutions vary rapidly {near} certain thresholds. Whether those jumps are smooth but fast varying, or are truly discontinuous, changes the way they are treated mathematically. It is natural to assume that a limit exists in which `fast-but-smooth' becomes truly `discontinuous', but how to characterise that limit is by no means {obvious}. 

Analytically this might not be surprising, because in many cases the limit may be singular (for an interesting aside on singular limits see \cite{berry94red}), meaning that solutions of the smooth system do not limit to the piecewise-smooth system in a regular way. Computationally the difficulty in studying the nonsmooth limit of smooth systems lies in the fact that, by their very nature, they become extremely stiff and {numerical methods fail to converge} in the limit of interest {\cite{guck00}}. As a result, the relation between smooth and piecewise-smooth dynamical systems theory remains poorly understood. The recent growth of piecewise-smooth dynamical systems theory has produced numerous forms of discontinuity-induced singularities, bifurcations, and chaos. It is interesting to ask whether these have any counterpart in smooth systems, where novel behaviours have also been attributed to abrupt change. For this purpose a relation between smooth and piecewise-smooth is necessary. (Examples of the discontinuity-induced phenomena of interest are sliding bifurcations \cite{bc08,bkn02,jh09}, explosions \cite{j11prl}, grazing singularities \cite{f88,t93}, and non-deterministic chaos \cite{cj09}. Examples of the smooth system phenomena of interest are canard explosions \cite{b81} and mixed-mode oscillations \cite{des12,krupa08,pop08}). 

In \cite{jd10}, ideas from singular perturbations, nonsmooth dynamics, and nonstandard analysis, were combined to develop a method called {\it pinching}, which approximates a fast-but-smooth change by a discontinuity. Pinching characterises the dynamics in the nonsmooth limit at least qualitatively, and, as we show here, even quantitatively, by preserving certain {singularities} and associated geometry. The present paper investigates the method by applying it to a known system with slow-fast dynamics, characterized by {an invariant manifold of slow dynamics} which loses {normal} hyperbolicity, taking the form of the so-called {\it folded-node} system in a singular limit. A relation between so-called {\it canard} phenomena in smooth and piecewise-smooth systems is derived. 

Pinching is a way of deriving piecewise-smooth models that capture key geometry of singular perturbation problems. While the method so-far developed does not provide a rigorous approximation {in the analytic sense}, it faithfully captures singularities and bifurcations that arise when a smooth system suffers rapid change, encapsulated in a piecewise-smooth differential equation, capable of providing accurate estimates of bifurcation parameters, as shown previously in a study of the van der Pol system in \cite{jd10}. 

In the remainder of this section we introduce the {canonical model used to study the folded node, a set of ordinary differential equations with} a singular perturbation parameter $\eps$. In \sref{sec:sans} we apply pinching to obtain a piecewise-smooth system dependent on $\eps\neq0$, whose phase portrait resembles the smooth system's {singular} limit $\eps=0$. Neither the pinched system, nor the singular limit of the smooth system, accurately represent the dynamics found when $\eps$ is small but nonzero. In \sref{sec:1st} we motivate pinching better by first making an exponential rescaling of the phase space, yielding a system with similar qualitative features to \sref{sec:sans}{; we then give numerical evidence, by means of boundary-value problem continuation, of the continuum of canards emerging in the nonsmooth limit}. In \sref{sec:2nd} we show how this can be improved, by shifting the focus of the exponential rescaling, and pinching again. This second approximation possesses more intricate dynamics, which we show is in one-to-one correspondence with solutions of the smooth system for $\eps\neq0$. We make some closing remarks in \sref{sec:conc}, including a simple smoothing of the pinched system from \sref{sec:sans} that shows how the number of canards in a flow increases, through a series of bifurcations, as a change at some threshold tends towards discontinuity. 

\begin{figure}[h!]\centering\includegraphics[width=0.4\textwidth]{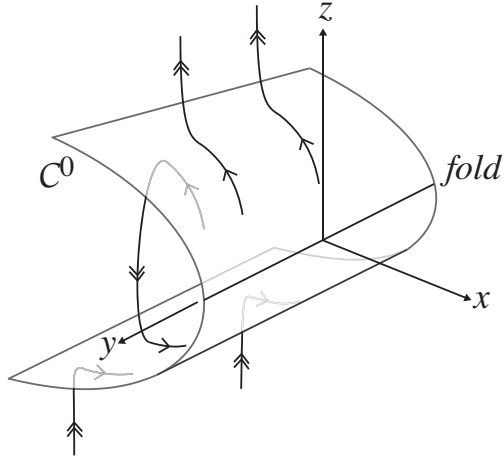}
\vspace{-0.3cm}\caption{\footnotesize\sf Critical manifold {${\cal C}^0$}, and dynamics away from the fold (the $y$-axis). Solutions evolve {quasi-}vertically on the fast timescale (double arrows) to/from the $\eps$-neighbourhood of the critical manifold, where the slow timescale dominates (single arrows). }\label{fig:fold}\end{figure}

We begin with a slow-fast system {whose key features} can be considered fundamental to the understanding of the canard phenomenon. {The system is a set of ordinary differential equations for two slow variables, $x$ and $y$, and a fast variable, $z$. (Some features of these systems can be generalised to arbitrarily many fast and slow variables, see \cite{brons06,w12}).} The timescale separation is introduced by a parameter $\eps$ satisfying $0<\eps\ll1$. A general such system can be written
\eqblock{ode}{
\dot x&=&g_1(x,y,z;\eps)\;,\\
\dot y&=&g_2(x,y,z;\eps)\;,\\
\eps\dot z&=&\;h(x,y,z;\eps)\;,
}
in terms of smooth functions $h,g_1,$ and $g_2$, where the dot denotes differentiation with respect to the (slow) time $t$. For $|h|>\eps$, it is easy to see that $\dot z$ is much larger than $\dot x$ and $\dot y$, so on the fast timescale solutions are attracted to, or repelled from, an $\eps$-neighbourhood of the  null surface of the fast variable, 
\begin{equation}\label{C0}
{\cal C}^0:=\{h(x,y,z,\eps)=0\}\;,
\end{equation}
known as the {\it critical manifold}; see \fref{fig:fold}. As shown by Fenichel \cite{f79}, within that neighbourhood {lie (generally non-unique)} perturbation{s} of the critical manifold, called the {\it slow manifold{s}}, hypersurface{s} which {are (locally)} invariant in the flow, and {are} hyperbolically attracting or repelling provided that $\partial h/\partial z\neq0$. {The notion of local invariance means that solutions can leave a slow manifold {\it only} at its boundary, if one exists \cite{j95}.} 

{Normal h}yperbolicity of the {critical} manifold is lost if $\partial h/\partial z$ vanishes. This scenario occurs generically in three dimensions when there is a fold {in} {${\cal C}^0$} with respect to the flow (see \fref{fig:fold}). A fold is a set of points where the conditions
\begin{eqnarray}
h=\frac{\partial\;}{\partial z}h=0,\label{fold}\\
{\left\{\frac{\partial\;}{\partial x}h,\frac{\partial\;}{\partial y}h\right\}}\neq0\neq \frac{\partial^2\;}{\partial z^2}h\label{foldndg},
\end{eqnarray}
are satisfied, with the inequalities ensuring that the critical points along $h=0$ are indeed folds, and not higher order degeneracies. The flow's projection {o}nto the $x$-$y$ plane is transverse to the fold except at points where
\begin{equation}\label{foldsing}
h=(g_1,g_2)\cdot\bb{\frac{\partial\;}{\partial x},\frac{\partial\;}{\partial y}}h=0\;.
\end{equation}
The three conditions in \eref{fold} and \eref{foldsing}, subject to non-degeneracy conditions (\ref{foldndg}),  define an isolated point, which we can place at the origin $(x,y,z)=(0,0,0)$ for $\eps=0$, where the system can be transformed into the local canonical form \cite{ben90,sw01}
\eqblock{odefold}{
\dot x&=&by+cz\;+\;h.o.t.\\
\dot y&=&1\;+\;h.o.t.\\
\eps\dot z&=&x+z^2\;+\;h.o.t.\\
}
This provides a local model for the system of interest, and from now on we omit higher order terms ($h.o.t.$). In doing so we must take care, for we can scale $\eps$ out of the leading order terms by sending $\cc{x,y,z,t}\mapsto\cc{x\eps,y\sqrt\eps,z\sqrt\eps,t\sqrt\eps}$ (following \cite{w05}), or $\cc{x,y,z,b,c}\mapsto\cc{x\eps^2,y,z\eps,b\eps^2,c\eps}$, and while the result is a formally correct leading order system, the $\eps$-orders of the new variables $(x,y,z,t)$ must be considered in any subsequent approximation. Since the order of certain quantities is essential here, we do not make such a substitution at this stage. 

\subsection{Projection onto the critical manifold}\label{sec:h=0}

This paper presents a method for characterising slow-fast dynamics for a system with small $\eps$, without taking $\eps$ to zero. The more common method of studying a system such as (\ref{odefold}) involves first taking the singular limit, $\eps=0$, and subsequently considering $\eps>0$ as a perturbation. Certain features of the singular limit will appear in a different guise (for $\eps\neq0$ in fact) in our nonsmooth approach later, so let us review these first. 

\begin{figure}[h!]\centering\includegraphics[width=0.9\textwidth]{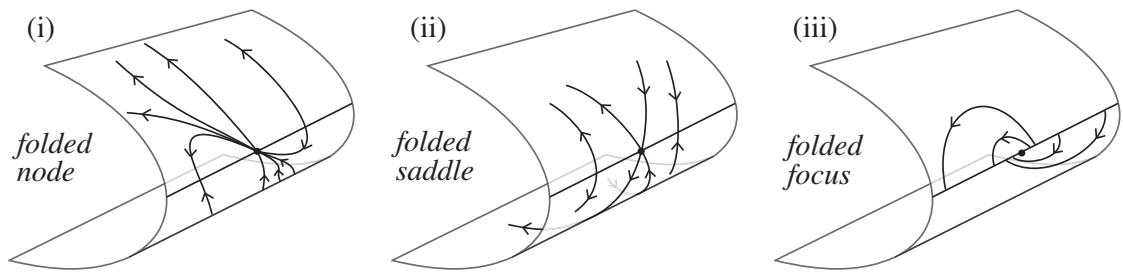}
\vspace{-0.3cm}\caption{\footnotesize\sf Projected onto the critical manifold, the flow is a singular scaling of (or a {\it folded}) (i) node, (ii) saddle, or (iii) focus.}\label{fig:foldsing}\end{figure}

Setting $\eps=0$ (and neglecting higher order terms) reduces \esref{odefold} to a differential-algebraic equation 
\begin{eqnarray*}
\dot x&=&by+cz\;,\\
\dot y&=&1\;,\\
0&=&x+z^2\;,
\end{eqnarray*}
called the {\it reduced system}, see e.g. \cite{ben90,sw01}. Its solutions are restricted to the critical manifold ${\cal C}^0$, on which $h=x+z^2=0$. This flow therefore satisfies $\dot h=0$ which, using \esref{odefold}, gives
$$0=\dot h=by+cz+2z\dot z\quad\Rightarrow\quad\dot z=-\frac{by+cz}{2z}\;.$$
Combining this with the second component $\dot y=1$ from \esref{odefold} yields a dynamical system on ${\cal C}^0$, given by
\begin{equation}\label{slowsub}
\bb{\begin{array}{c}\dot y\\\dot z\end{array}}=\frac{-1\;}{2z}\bb{\begin{array}{cc}0&-2\\b&c\end{array}}\bb{\begin{array}{c}y\\z\end{array}}\;.
\end{equation}
This system is undefined at $y=z=0$ but, importantly, arbitrarily close to this point the righthand side of \eref{slowsub} is generally nonzero. The local phase portrait is determined by the $2\times2$ in \eref{slowsub}, whose trace and determinant are respectively
\begin{equation}\label{dsinglambda}
\lambda_1\lambda_2=2b\qquad\mbox{and}\qquad\lambda_1+\lambda_2=c\;,
\end{equation}
where $\lambda_{1,2}$ are solutions of the characteristic equation $0=\lambda^2-c\lambda+2b$ assigned such that $|\lambda_1|\ge|\lambda_2|$. Let $\mu=\lambda_2/\lambda_1$. The singularity is then classified (see for example \cite{sw01}) as:
\begin{eqnarray}
\mbox{a folded node }&\mbox{if}&\mu>0\;,\qquad(\;\lambda_1\lambda_2>0\;)\\
\mbox{a folded saddle }&\mbox{if}&\mu<0\;,\qquad(\;\lambda_1\lambda_2<0\;)\\
\mbox{a folded focus }&\mbox{if}&\mu\in\mathbb C\;.\qquad(\;\lambda_1=\lambda_2^*\;\;\;)
\end{eqnarray}
The three cases are illustrated in \fref{fig:foldsing}. Their names reflect the fact that, if we omit the prefactor $-1/2z$ from \eref{slowsub}, the remaining linear system 
$$\bb{\begin{array}{c}\dot y\\\dot z\end{array}}=\bb{\begin{array}{cc}0&-2\\b&c\end{array}}\bb{\begin{array}{c}y\\z\end{array}}$$
has an equilibrium at the origin, and this is a node if $\mu>0$, a saddle if $\mu<0$, or a focus if $\mu\in\mathbb C$. The term `folded' is required because this linear system is topologically, but not dynamically, equivalent to \eref{slowsub}, being obtained from it by a time scaling $t\mapsto-t/2z$. This scaling changes sign with $z$, and is singular at $z=0$. As a result, the dynamics is similar in $z<0$ on the attracting branch of the critical manifold ${\cal C}^0$, and is similar up to time-reversal in $z>0$ on the repelling branch of ${\cal C}^0$. The singularity of the scaling at $z=0$ means that the equilibrium of the linear system is not an equilibrium of \eref{slowsub}, instead referred to as a {\it folded equilibrium}, which, unlike an equilibrium, solutions can cross through in finite time. Solutions of \eref{slowsub} that cross through $y=z=0$ are shown in \fref{fig:foldsing}, with infinitely many in (i), only two in (ii), and none in (iii). Those that pass from $z<0$ to $z>0$ through the folded node or folded saddle correspond to canards, which are discussed in the next section.

The classification into folded node/saddle/focus, for different values of $\mu$, is also used to classify the \esref{odefold} for $\eps$ nonzero. 
The remainder of this paper will be concerned solely with the folded node case, $\mu>0$. Substituting $\mu=\lambda_2/\lambda_1$ and \eref{dsinglambda} into \eref{odefold}, then scaling {$(x,y,z,t,\eps)\mapsto ( x/\lambda_1^2,y/\lambda_1^2 ,- z/\lambda_1, t/\lambda_1^2,-\eps/\lambda_1^3)$}, gives to leading order
\eqblock{muform}{
\dot x&=&\frac\mu2y-(1+\mu)z\;,\\
\dot y&=&1\;,\\
\eps\dot z&=&x+z^2\;.
}
Since much of the foregoing analysis will involve a natural length scale $|h|=\eps$, for convenience we define a new function
\begin{equation}\label{hbar}
 u (x,y,z;\eps)=h(x,y,z;\eps)/\eps\;,
\end{equation}
then in the variables $\cc{ u ,y,z}$, the \esref{muform} becomes
\eqblock{straight}{
\eps\dot  u &=&\frac\mu2y{-}(1+\mu)z+2z u \;,\\
\dot y&=&1\;,\\
\dot z&=& u \;.
}
The remainder of the paper is a study of this system for $\mu>0$ and $\eps>0$. 
Before we apply pinching to this system, we must state some {preliminaries} concerning the so-called {\it canard} type solutions that make it so interesting.

\subsection{Canards}\label{sec:can}

Throughout this paper we define a {\it canard} as a solution that evolves from an attracting invariant manifold, to a repelling invariant manifold, via some singularity that faciliataes the transition. In the traditional setting of a smooth system with slow and fast timescales, the invariant manifolds are surfaces of slow dynamics, as we introduce in \eref{cant} below. In the setting of discontinuous systems as we introduce in \sref{sec:sans}, the invariant manifolds are regions where solutions slide along the switching manifold (the discontintuity set). Strictly speaking, we use the term `canard' exclusively for {\it maximal canards}, which are the solutions of the types above that spend the maximum possible time on the repelling invariant manifold;in the setting described in this paper, this time can be regarded as infinite. 

For real $\mu$ and noting $|\mu|<1$ by definition, the ratios $z/y=\mu/2$ and $z/y=1/2$ are satisfied by the weak and strong eigendirections associated with the folded equilibrium of \eref{slowsub}. The solutions along these directions are canards, known as the weak and strong {\it singular canards} of the singular ($\eps=0$) system. 

The non-singular (i.e. $\eps\neq0$) system (\ref{straight}) also has two particular solutions satisfying ${z}/{y}=\mu/2$ and $z/y=1/2$. 
We label these $\gamma^{wk,st}=\cc{ u (t),\;y(t),\;z(t)}$, and solve \eref{straight} to find
\begin{equation}\label{cant}
\gamma^{wk}(t)=\cc{\frac\mu2,\;t,\;\frac{\mu}{2}t\;},\qquad
\gamma^{st}(t)=\cc{\frac12,\;t,\;\frac{1}{2}t\;}\;.
\end{equation}
These solutions are canards, since as $t\rightarrow\pm\infty$ they lie in an $\eps$-neighbourhood of the attracting and repelling branches of the critical manifold, implying that they tend towards attracting and repelling slow manifolds. They form simple curves that are $\eps$-close to the weak and strong canards of the singular system, therefore $\gamma^{wk}$ is called the {\it weak {primary} canard}, and $\gamma^{st}$ the {\it strong {primary} canard}, (though where possible without ambiguity we omit the word `primary'). 

\begin{figure}[h!]\centering\includegraphics[width=0.9\textwidth]{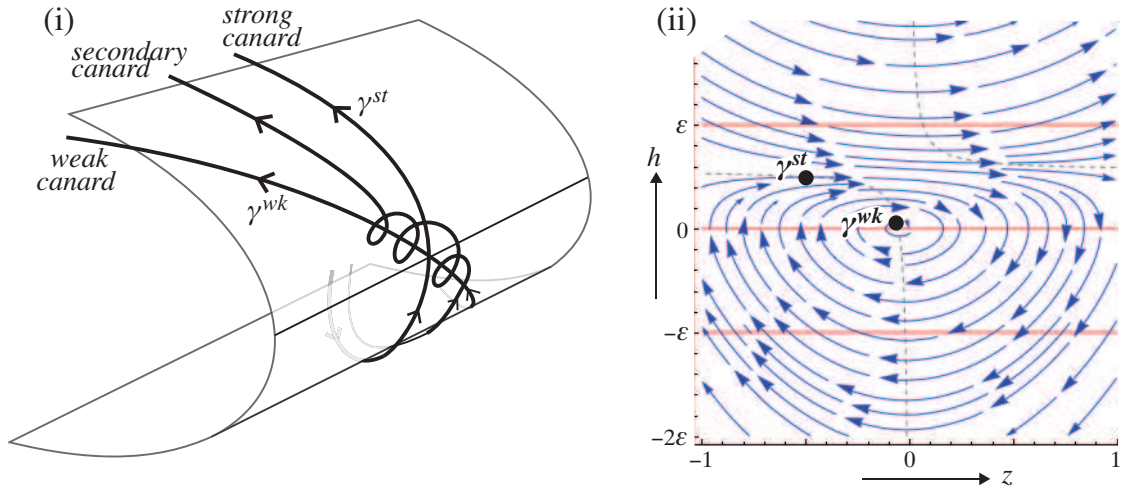}
\vspace{-0.3cm}\caption{\footnotesize\sf (i) Primary and secondary canards in the neighbourhood of the critical manifold. Secondary canards connect an $\eps$-neighbourhood of the attracting and repelling branches of the critical manifold, their tails lying parallel to the strong {(primary)} canard $\gamma^{st}$. Near the singularity they rotate around the weak (primary) canard $\gamma^{wk}$. (ii) The flow circulating around the weak canard: simulation of \eref{straight} in the plane $y=-1$.}\label{fig:canards}\end{figure}

It is clear from the local phase portrait (see \fref{fig:foldsing}(i)) that, besides the primary canards, the singular system contains a whole family of canard solutions through the singularity, forming a continuum between the weak and strong solutions. For the non-singular ($0<\eps\ll1$) system there may exist a number of other canards, termed {\it secondary canards}, though unlike the singular system they will generally be finite in number, and of a more complicated topology (see for example \cite{w05,des10}). Secondary canards have been shown \cite{w05,des10} to rotate around the weak {(primary)} canard near the origin, and to asymptotically align with the strong {(primary)} canard as $t\rightarrow\pm\infty$, as sketched in \fref{fig:canards}. They are neither easy to express in closed form, nor easy to simulate numerically. To study secondary canards analytically, Wechselberger \cite{w05} applies a parameter blow-up, then moves to cylindrical coordinates centred on the weak canard, and takes the variational equation along the weak canard to obtain a Weber equation, whose solutions {describe small oscillations that} the secondary canards perform around the weak canard. In this paper pinch the $\eps$-neighbourhood of the critical manifold onto the manifold itself, thereby deriving a piecewise-smooth system that exhibits analogous behaviour. 

For $\eps$ nonzero, a first approximation for the slow manifolds is that they lie in an $\eps$ neighbourhood of the critical manifold (at least where it is normally hyperbolic). A better approximation, and one we will use later, is to note that the slow dynamics lies not on, but is stationary with respect to, the critical manifold, and hence lies close to (actually in an $\eps^2$-neighbourhood of) the nullcline $\dot u=0$. Solving $\dot u=0$ in \eref{straight} gives the surface 
\begin{equation}\label{Py}
{\cal P}_y=\left\{( u ,y,z)\in\mathbb R^3:\; u =\frac{\mu+1}2-\frac{\mu y}{4z}\right\}\;.
\end{equation}
These two approximations ($\eps$-close to ${\cal C}^0$ and $\eps^2$-close to ${\cal P}_y$) must be consistent, that is, the surface ${\cal P}_y$ can only approximate a slow manifold where it lies in the $\eps$-neighbourhood of the critical manifold. Near $z\approx=0$ this cannot hold, while for large $z$ the surface ${\cal P}_y$ is approximated by 
\begin{equation}\label{P0}
{\cal P}_0=\cc{( u ,y,z)\in\mathbb R^3: u =\frac{\mu+1}2\;}\;,
\end{equation}
and for this to lie within $\eps$ of the critical manifold, which corresponds to the region $|u|<1$, we must clearly have $|(\mu+2)/2|<1$. 
Combining this with the folded node condition $\mu>0$ gives the range of permitted values for the parameter $\mu$ as
\begin{equation}\label{cancon2}
0<\mu<1.
\end{equation}
This motivates the choice we made earlier of defining $\mu$ as $\lambda_2/\lambda_1$ rather than its reciprocal. As an illustration of the condition above, note that the weak and strong canards are two particular solutions that lie on the slow manifolds. These lie on ${\cal P}_y$ because their distances $u$ from the critical manifold are fixed at $u=\mu/2$ and $u=1/2$, and moreover $|\mu|<1$ guarantees that both of these lie in the required neighbourhood, since $|\mu/2|<1$ and $|1/2|<1$. 

\section{Pinching}\label{sec:sans}

Pinching, at least in the form used here, was introduced in \cite{jd10}. To illustrate the method, we first demonstrate it in its crudest form. Essentially, we assume the system is dominated by fast dynamics for $| u |>1$, which we leave untouched. We assume that slow dynamics dominates in $| u |<1$, in such a way that dynamics in the $ u $ direction can be neglected, so we collapse the entire neighbourhood $| u |<1$ onto the critical manifold $ u =0$. 

The `pinch` is enacted by a piecewise-smooth transformation of the variables, in this case introducing a new variable
\begin{equation}\label{sans}
U= u -{\sf sign}( u )\;,
\end{equation}
for $| u |>1$, and omitting $| u |<1$. \Erefs{straight} then become a piecewise smooth system
\eqblock{pinchsans}{
\eps\dot{U}&=&\frac\mu2 y-(\mu+1)z+2z\;{\sf sign}(U)\;\;+\;\;\ord{zU}\;,\\
\dot y&=&1\;,\\
\dot z&=&U+{\sf sign}(U)\;.
}
where 
\begin{equation}\label{sgn}
{\sf sign}(U)\in\left\{\begin{array}{lll}U/|U|&{\rm if}&U\neq0,\\\sq{-1,{+}1}&{\rm if}&U=0.\end{array}\right.
\end{equation}
The dynamical theory of such differential inclusions was described by Filippov \cite{f88}. For $U\neq0$, \eref{pinchsans} specifies the fast dynamics uniquely. On $U=0$, called in piecewise-smooth dynamics the {\it switching manifold}, the righthand side is set-valued. The dynamics it gives rise to, however, is rather simple to describe. 

\begin{figure}[h!]\centering\includegraphics[width=0.9\textwidth]{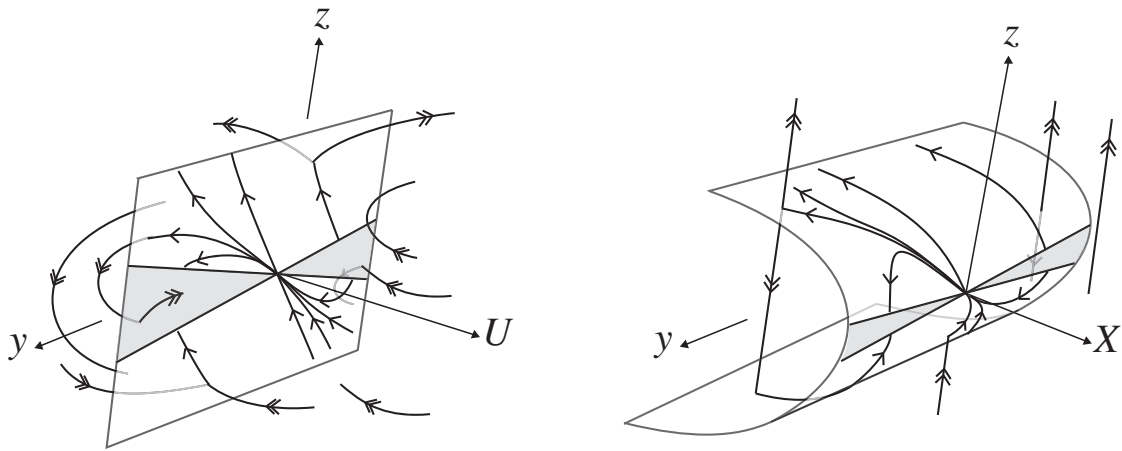}
\vspace{-0.3cm}\caption{\footnotesize\sf Pinched system, in $(U,y,z)$ coordinates (left) and corresponding $(X,y,z)$ coordinates (right) if $X=x-{\sf sign}(x)$. The upper and lower regions are the repelling and attracting sliding regions, separated by a `bow-tie'  (shaded) made up of crossing regions. Double arrows indicate trajectories outside the surface, single arrows indicate sliding trajectories.}\label{fig:sans}\end{figure}

Consider a point $p$ on the switching manifold, so $U|_p=0$. If $\dot U$ is nonzero there and its sign does not change with the sign of $U$, then $p$ is the arrival point of a solution lying on one side of $U=0$, and the departure point of a solution lying on the other side. Concatenating the two solutions gives a unique, continuous but non-differentiable, solution, that crosses the switching manifold at $p$. This therefore takes place on $U=0$ where $4z^2<\bb{\frac\mu2y-(\mu+1)z}^2$ (the `bow-tie' region in \fref{fig:sans}). 

If $\dot{U}$ changes sign at $p$ then the flow cannot cross through the switching manifold there. Two solutions of \eref{pinchsans} meet at $p$, each arriving from either side of $U=0$ (the lower region on $U=0$ in \fref{fig:sans}), or each departing (the upper region on  $U=0$ in \fref{fig:sans}). They occupy the regions on $U=0$ where $4z^2>\bb{\frac\mu2y-(\mu+1)z}^2$; the region where $z<0$ and $4z^2-\bb{\frac\mu2y-(\mu+1)z}^2>0$ attracts the flow outside the switching manifold, while the region where $z>0$ and $4z^2-\bb{\frac\mu2y-(\mu+1)z}^2>0$ repels it. These are the pinched analogues of a point $p$ on the attracting and repelling slow manifolds, respectively, of \esref{straight}. 
To find the flow on them we solve as we did for the slow flow projected onto $ u =0$ in \sref{sec:h=0}, fixing $0=\dot{U}=\frac\mu2 y-(\mu+1)z+2z\dot z$. Solving for $\dot z$ gives dynamics on $U=0$ defined by
\begin{equation}\label{musliding}
\bb{\begin{array}{c}\dot y\\\dot z\end{array}}=\frac{-1\;}{2z}\bb{\begin{array}{cc}0&-2\\\frac\mu2&-(\mu+1)\end{array}}\bb{\begin{array}{c}y\\z\end{array}}\;.
\end{equation}
We say that \esref{musliding} defines {\it sliding dynamics} on $U=0$, and its solutions are known as sliding orbits. 
This is clearly sensible as an analogue of the slow dynamics, since it is equal to \eref{slowsub}, which describes the original system projected onto the critical manifold. 

Finally, notice what happens at the boundaries of the regions of sliding and crossing, where $U=0$ and $\frac\mu2 y-(\mu+1)z\pm2z=0$. The field in $U>0$ is tangent to the switching manifold along $\frac\mu2 y-(\mu+1)z+2z=0$, and curves away from the manifold since it satisfies $\ddot U=\frac\mu2-(\mu+1)+2=1-\frac\mu2>0$. Thus the flow in $U>0$ at such points carries a solution away from the switching manfiold. The field in $U<0$ is tangent to the switching manifold along $\frac\mu2 y-(\mu+1)z-2z=0$, and curves towards the manifold since $\ddot U=\frac\mu2+(\mu+1)+2=3(\frac\mu2+1)>0$, so solutions can only enter the switching manifold at such points.

All this gives the simple dynamical portrait shown in \fref{fig:sans}. A feature of this system is that it contains solutions corresponding to the weak and strong canards, and a continuum of other canards between them. Evidently, recalling \sref{sec:can}, the sliding dynamics corresponds to represent the continuum of canards that exist in the smooth system's singular limit $\eps=0$. But the pinched system also allows crossing of the switching manifold, which the singular ($\eps=0$) smooth system does not, so it can only apply to $\eps\neq0$. 
It will emerge that the correct canard structure can be captured by modifying the pinch slightly (looking ahead, we do this by pinching around the surface (\ref{P0}) instead of the critical manifold). In \sref{sec:1st} we re-consider the choice of coordinates above, inserting a preparatory step (a ``microscope") that more fully motivates the process of pinching, before changing the focus of the microscope and the pinch to capture more precisely the secondary canard structure in \sref{sec:2nd}.

\section{First approximation: a continuum of canards}\label{sec:1st}

\begin{figure}[h!]\centering\includegraphics[width=0.9\textwidth]{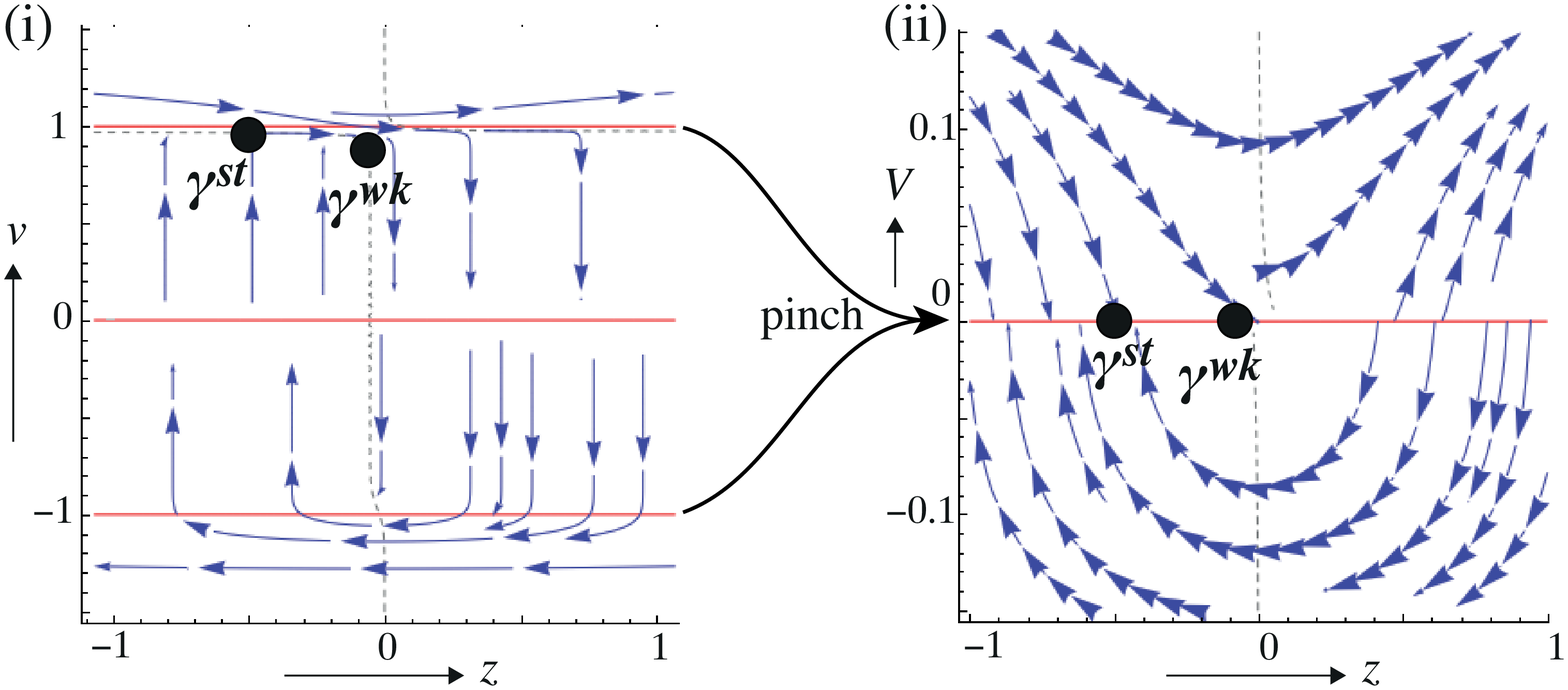}
\vspace{-0.3cm}\caption{\footnotesize\sf  The microscope and pinch. (i) Flow in the microscope \esref{zoom1} simulated in the plane $y=-1$ with $\eps=0.05$ and $\mu=1/8.5$. (ii) Flow in the \esref{pinch1} obtained by pinching together the surfaces $v=\pm1$ in (i), with the vertical axis rescaled for clarity. The primary canards $\gamma^{st,wk}$ are indicated.}\label{fig:0}\end{figure}

The microscope, introduced in the context of nonstandard analysis \cite{b81}, is an exponential scaling of variables that attempts to resolve the interaction between slow and fast dynamics. We use it here to motivate the pinch that will follow. 
Let
\begin{equation}
v= u ^{[\eps]}\;:=\;| u |^\eps {\sf sign}\; u \;,
\end{equation}
in terms of which \eref{straight} becomes
\eqblock{zoom1}{
\dot v&=&v\cc{2z+\bb{\frac\mu2y-(\mu+1)z}v^{-{1}/{[\eps]}}}\;,\\
\dot y&=&1\;,\\
\dot z&=&v^{{1}/{[\eps]}}\;.
}
This is simulated in \fref{fig:0}(i). The weak and strong canards now satisfy $v^{1/[\eps]}=z/y=\mu/2$ and $v^{1/[\eps]}=z/y=1/2$ respectively (by direct calculation). 
The projection onto $v=0$ gives the same 
$y$-$z$ system (\ref{musliding}) (or equivalently (\ref{slowsub})) as before the microscope. 

For $|v|<1$ and $\eps\ll1$, the flow of \esref{zoom1} is dominated by the term  
$v^{1-\frac1{[\eps]}}$ in the equation for $\dot v$, and so lies close to a set of fibres with constant $y$ and $z$, connecting the surfaces $v=\pm1$. To approximate this, we pinch the two surfaces $v=\pm1$ together using the transformation
\begin{equation}\label{vsans}
V=v-{\sf sign}\;v
\end{equation}
for $|v|>1$, giving the system
\eqblock{pinch1}{
\dot V&=&(V+{\sf sign}V)\cc{2z+\bb{\frac\mu2y-(\mu+1)z}(V+{\sf sign}V)^{-{1}/{[\eps]}}}\;,\\
\dot y&=&1\;,\\
\dot z&=&(V+{\sf sign}V)^{{1}/{[\eps]}}\;,}
as simulated in \fref{fig:0}(ii). 
Note that prior to pinching, both the weak and strong canards lie inside the region $|v|<1$, therefore they are not part of the \esref{pinch1} for $V\neq0$. Instead, they now lie on the switching manifold $V=0$. At the switching manifold, $V=0$, \eref{pinch1} reduces to
\eqblock{pinch1V0}{
\dot V&=&\frac\mu2y-(\mu+1)z+2z\;{\sf sign}V\;,\\
\dot y&=&1\;,\\
\dot z&=&{\sf sign}V\;,}
which is equivalent to \esref{pinchsans} up to an $\eps$ scaling in $y,z,t$. The crossing and sliding dynamics on the switching manifold are therefore exactly as described in \eref{musliding} for \esref{pinchsans}. The microscope, given by the transformation $ u \mapsto v$, provides the motivation for the pinch $ u \mapsto V$. The result is a system of `fast' dynamics for $V\neq0$ given by \eref{pinch1}, and slow dynamics that either crosses the sliding manifold, or slides along it as described in \eref{musliding} for \esref{pinchsans}. 
An immediate consequence is that, similar to the previous section (where the pinch was applied without a microscope), pinching here gives a continuum of canards solutions. In the next section the microscope is used more powerfully, to resolve the different canards. 

\begin{figure}[h!]\centering\includegraphics[width=0.4\textwidth]{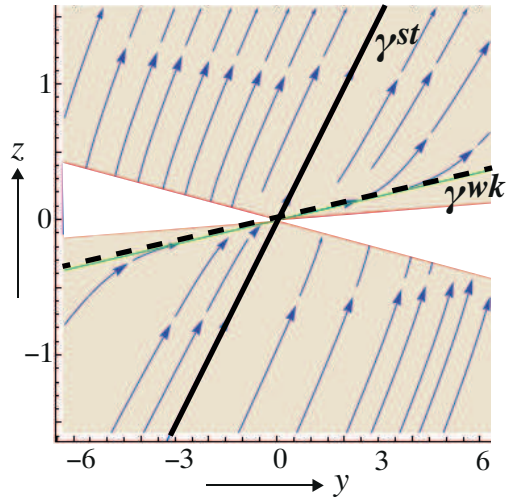}
\vspace{-0.3cm}\caption{\footnotesize\sf Sliding flow after the pinch, given by \eref{musliding} in the regions $4z^2>\bb{\frac\mu2y-(\mu+1)z}^2$ on $V=0$. Showing the strong (bold) and weak (dashed) primary canards. Orbits to the right of the strong canard in the upper half space (the repelling sliding region), and to the left of the strong canard in the lower half space (the attracting sliding region), form a continuum of canards. }\label{fig:fil0}\end{figure}

\subsection{Computations in the nonsmooth limit}\label{sec:numerics}

Before resolving the secondary canards by adjusting the pinch, one may immediately ask whether the original canard structure is restored when we smooth out the piecewise-smooth model (\ref{pinchsans}). The obvious way to achieve this is to replace the ${\sf sign}$ function with a  sigmoid function. Neglecting the $\ord{zU}$ term in \eref{pinchsans}, we directly replace $\sf sign$ with 
$${\sf sign}(U)\mapsto\tanh(kU)\;,$$
for large positive constant $k$, and {analyse} the resulting system {numerically} using a technique {based on the numerical continuation of parametrised families of two-point boundary-value problems} as {presented} in \cite{des10}. The piecewise-smooth \esref{pinchsans} is obtained in the limit $k\rightarrow\infty$, and the original smooth \esref{straight} is regained as $k\rightarrow0$.

{The computations below are all made for $\mu=1/8.5$. \Fref{fig:numerics1} shows attracting and repelling slow manifolds for different values of $k$, found by computing solutions that pass between lines chosen on the attracting and repelling branches of the critical manifold. Canards occur where the attracting and repelling branches intersect transversally. The magnified images show an increasing number of rotations around the weak canard, and an increasing number of intersections (i.e. canards), with increasing $k$, that is, as the smoothing function $\tanh(kU)$ approaches ${\sf sign}(U)$. 

\begin{figure}[h!]\centering\includegraphics[width=\textwidth]{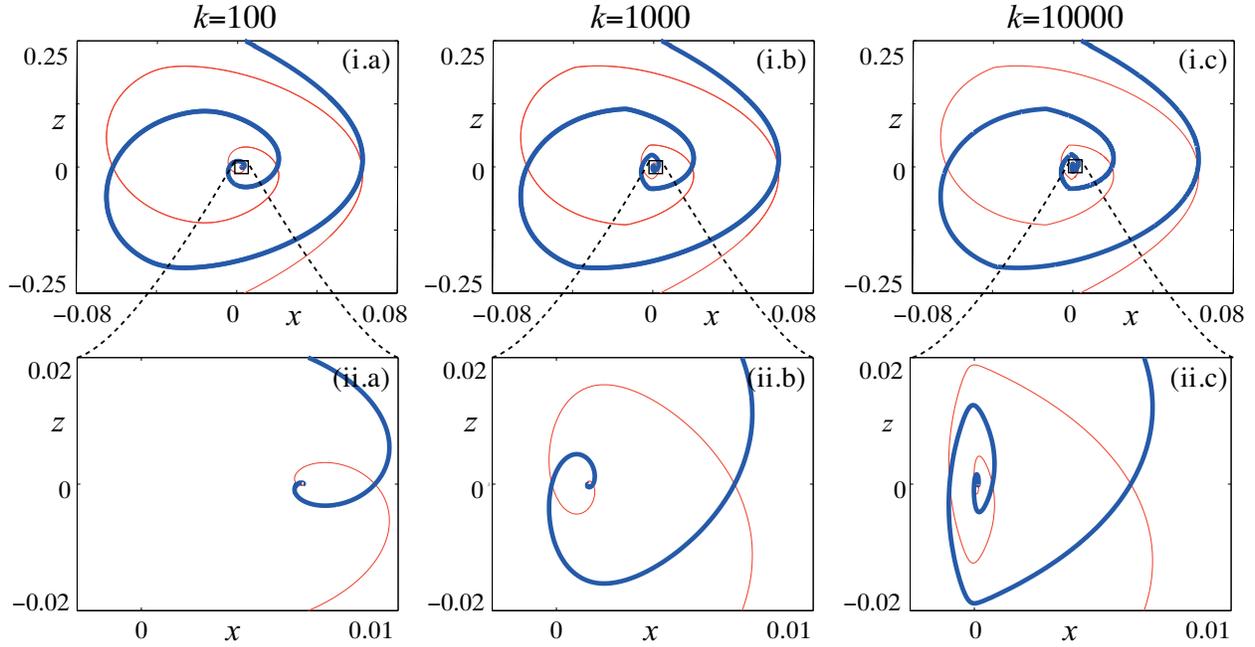}
\vspace{-0.3cm}\caption{\footnotesize\sf {Slow manifolds} in a regularization of the piecewise-smooth \esref{pinchsans},  plotted in the section {$\{y=0\}$} {for different values of the stiffness parameter $k$}. (i) Computations of the slow manifolds, and (ii) magnifications showing the increasing numbers of intersections with increasing $k$. Each intersection of the slow manifolds corresponds to a canard. (Colour online indicates the repelling (thicker/blue) and attracting (thinner/red) slow manifolds).}\label{fig:numerics1}\end{figure}

{\Fref{fig:numerics2} shows ten branches of canards {($\gamma^1$ to $\gamma^9$ and $\gamma^{st}$) continued as the stiffness parameter $k$ varies (see \cite{des10} for details on such continuation for canards)}. Each point along the ten curves corresponds to a canard, formed by transversal intersection of the the repelling and attracting branches of slow manifolds from \fref{fig:numerics1}}. The curves give the maximum value of $x$ reached along each canard solution (compare for example to theoretical results in \cite{w05}). The strong canard is labeled {$\gamma^{st}$}, and nearby, we see the development of secondary canards, labelled $\gamma^1$ to $\gamma^9$, as $k$ increases from zero.

{The number of canards grows very quickly with $k$ at smaller values ($k\lesssim1000$). Each new branch of secondary canards emerges from one particular curve, labeled as the branch of weak canards, $\gamma^{wk}$. This branch could not be computed for all values of $k$, and is partly derived from the envelope of endpoints of the secondary canards, however its identification as the weak canard is supported for several reasons. Primarily, only the weak canard should coexist with $\gamma^{st}$ for all parameters, and furthermore the bifurcation of secondary canard branches fits with previous results for canards in the case of a folded node; we refer the reader to \cite{w05} for theoretical results and a sketch of the expected bifurcation diagram, and to \cite{d08} for a computed bifurcation diagram.} 

\begin{figure}[t!]\centering\includegraphics[width=\textwidth]{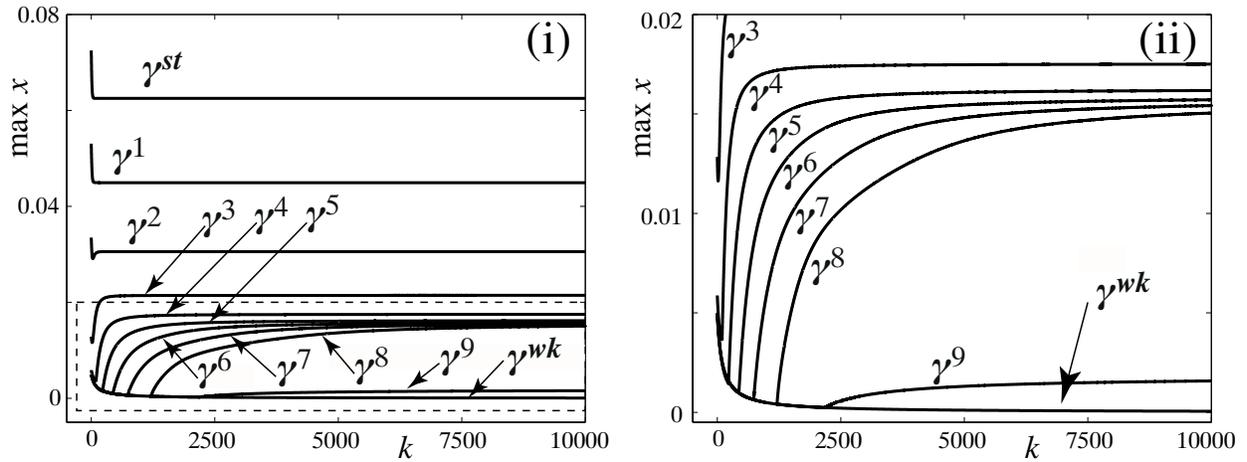}
\vspace{-0.3cm}\caption{\footnotesize\sf {Branches of canards in a regularization of the piecewise-smooth \esref{pinchsans}. For valid canard solutions, the maximum value of $x$ reached along the canard is plotted against the stiffness parameter $k$. The strong canard ($\gamma^{st}$), weak canard ($\gamma^{wk}$), and several secondary canards ($\gamma^1$-$\gamma^9$) are identified.}}\label{fig:numerics2}\end{figure}
Although only the first nine secondary canards are identified here, these computations suggest a trend towards infinitely many canards in the nonsmooth ($k\rightarrow\infty$) limit. This scenario fits with the model (\ref{pinch1}) obtained by pinching, suggesting that the piecewise-smooth \esref{pinch1} approaches the smooth \esref{straight}, in the limit where small $\eps$ leads to an increasingly sharp jump in dynamics near the critical manifold. {Even allowing for canards beyond those counted here, \fref{fig:numerics2} is expected to be incomplete, since for large $k$ ($k\gtrsim1000$) additional canard branches were detected in computations, born through folds instead of from the branch of weak canards, whose identity it currently unclear.} One can ask how these various features depend on the form of the sigmoid function that replaces ${\sf sign}(U)$, and whether they are represented in any way in the piecewise-smooth system. These are interesting problems for further study, with particular relevence to the study of uniqueness in regularization of piecewise-smooth systems. For the present paper, we now return to developing the pinch approximation.

\section{Second approximation: rotating canards}\label{sec:2nd}

The exponential scaling in the previous section attempts to resolve the slow dynamics in the neighbourhood of the critical manifold $ u =0$. However, after the microscope on $ u =0$, the flow is still seen to evolve fast towards another surface, the nullcline $\dot u =0$. An improvement on this method is therefore to take a new microscope and pinch, centred on this new surface. 

The nullcline where $\dot  u =0$ is given by ${\cal P}_y$ in \eref{Py}. Since this is not well defined at $z=0$, we approximate it by ${\cal P}_0$ given by \eref{P0}, and take a microscope on ${\cal P}_0$, by introducing a new variable
\begin{equation}\label{microscope2}
w=\bb{\frac u \eps-\frac{1+\mu}{2\eps}}^{[\eps]}\;,
\end{equation}
in terms of which the dynamical \esref{straight} becomes
\eqblock{zoom2}{
\dot w&=&w\cc{{2z}+\frac{\mu y}{2\eps}w^{-{1}/{[\eps]}}}\;,\\
\dot y&=&1\;,\\
\dot z&=& \eps w^{{1}/{[\eps]}}+\frac{1+\mu}2\;.
}
This is simulated in \fref{fig:1}(i). 
The weak and strong canards now satisfy $\cc{w,z}=\cc{-1/(2\eps)^\eps,\mu y/2}$ and $\cc{w,z}=\cc{-(\mu/2\eps)^\eps,y/2}$, respectively. 

The nullcline $\dot w=0$ is the curve
\begin{equation}
{\cal Q}_y=\cc{(w,z)\in\mathbb R^2:\; w=-\bb{\frac{\mu y}{4z\eps}}^{[\eps]}}\;.
\end{equation}
The flow projection onto $w=0$ has yet again the same $y$-$z$ system as before, namely \eref{musliding}. However, the fact that the flow organises around the nullcline ${\cal Q}_y$ as {apposed} to ${\cal P}_y$ gives different dynamics, as will be revealed by the pinch. 


Slow dynamics in \fref{fig:1}(i) is seen {numerically} to dominate in an $\eps^2$ neighbourhood of the nullcline $\dot h=0$, which corresponds to $| u -\frac{\mu+1}{2}|<\eps$ or $|w|<1$. In the first approximation of \sref{sec:1st}, slow dynamics is observed to dominate in an $\eps$-neighbourhood of the critical manifold, $|h|<\eps$, which corresponds to $| u |<1$ or $|v|<1$. An analytic explanation of why the slow neighbourhood is of order $\eps$ around the critical manifold, and $\eps^2$ around its associated nullcline, is outside the scope of the current paper, but deserves attention in future work. 

The pinch is now enacted similarly to the previous section, by introducing a new variable
$W=w-{\sf sign}\;w$, which gives
\eqblock{pinch2}{
\dot W&=& (W+{\sf sign}W)\cc{{2z}+\frac{\mu y}{2\eps}(W+{\sf sign}W)^{-{1}/{[\eps]}}}\;,\\
\dot y&=&1\;,\\
\dot z&=&\eps(W+{\sf sign}W)^{{1}/{[\eps]}}+\frac{1+\mu}{2}\;,}
as shown in \fref{fig:1}(ii).

\begin{figure}[h!]\centering\includegraphics[width=0.9\textwidth]{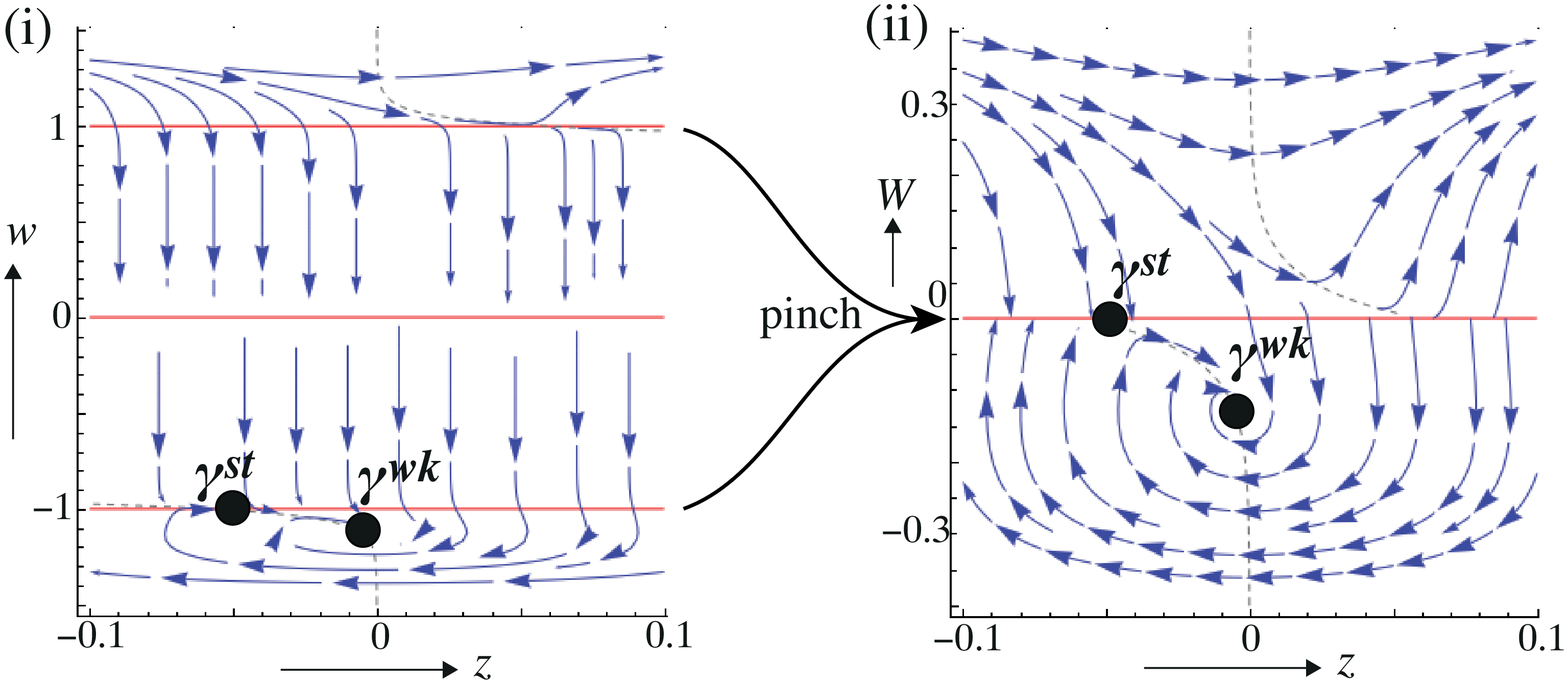}
\vspace{-0.3cm}\caption{\footnotesize\sf  The second microscope and pinch. (i) Flow in the microscope \esref{zoom2} simulated in the plane $y=-1$ with $\eps=0.05$ and $\mu=1/8.5$. (ii) Flow in the \esref{pinch2} obtained by pinching together the surfaces $w=\pm1$ in (i), with the vertical axis rescaled for clarity. The primary canards $\gamma^{st,wk}$ are indicated.}\label{fig:1}\end{figure}

The sliding flow is superficially given by the usual \eref{musliding}, as is found by solving for $\dot W=0$ on $W=0$. 
Crucially, however, we must consider the arrangement of the tangencies $\dot W=0$ on $W=0$, which give the boundaries of the sliding regions. The $W>0$ and $W<0$ subsystems in \eref{pinch2} are tangent to the switching manifold where $\dot W$ reaches zero as $W$ approaches zero from above or below. These tangencies lie along $\frac zy=-\frac\mu{4\eps}{\sf sign}W$, hence the sliding regions are found to be given by 
\begin{equation}\label{Wslidingregion}
|z/y|>\mu/4\eps\quad{\rm on}\quad W=0\;,
\end{equation} 
illustrated in \fref{fig:fil1} for different values of $\mu$. 

\begin{figure}[h!]\centering\includegraphics[width=\textwidth]{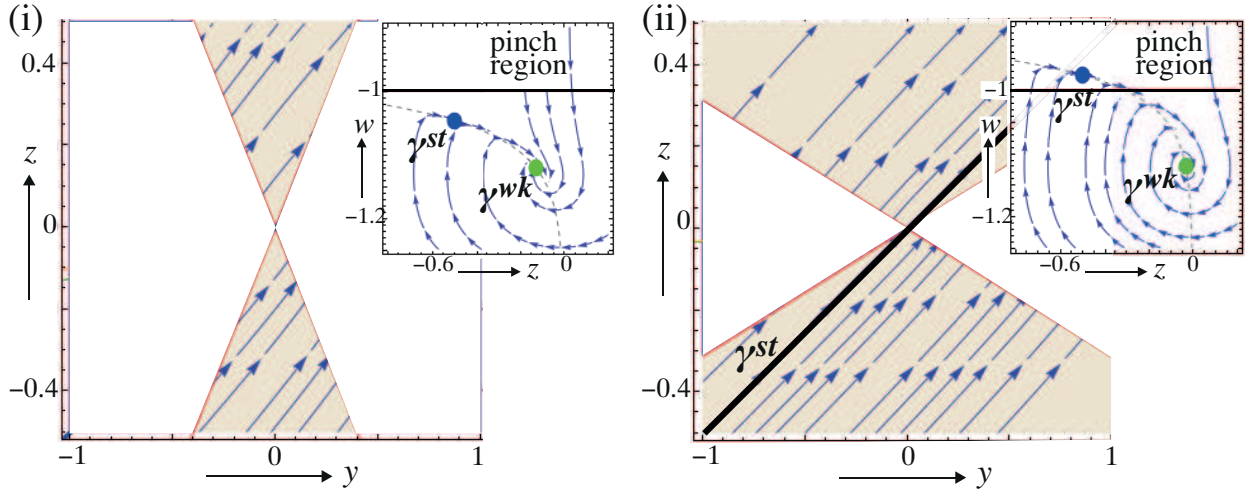}
\vspace{-0.3cm}\caption{\footnotesize\sf  Sliding flow after the second pinch, given by \eref{musliding} in the regions  $|z/y|>\mu/4\eps$ on $W=0$, for $\eps=0.05$. In (i) there are no canards for $\mu=1/4$, and in (ii) the strong canard can be seen in the sliding flow for $\mu=1/16$. Inset: the strong canard in the microscope system falls inside the pinch region in (ii) but not in (i).}\label{fig:fil1}\end{figure}

The curvature of the flow is specified by the second derivative, 
$\ddot W=
(1\pm2\eps)(\mu\pm2\eps)/2\eps$. Since we have $0<\mu<1$ and $0<\eps\ll1$, the flow in $W>0$ satisfies $\ddot W=(1+2\eps)(\mu+2\eps)/2\eps>0$ on its tangency line $z/y=\mu/4\eps$, and hence curves away from the switching manifold. The flow in $W<0$ satisfies $\ddot W=(1-2\eps)(\mu-2\eps)/2\eps$ on its tangency line $z/y=-\mu/4\eps$, and hence curves away from the switching manifold if $\mu<2\eps$, and towards it if $\mu>2\eps$. Although we are interested in arbitrarily small $\eps$, either of these can be satisfied for small enough $\mu$. 

This happens because the weak eigendirection of \eref{musliding} lies outside the sliding region. The strong eigendirection lies inside the sliding regions if $\mu>2\eps$, or outside the sliding regions if $\mu<2\eps$, meaning that in the former case the sliding dynamics captures no canards, and in the latter captures a single canard. (This is immediately in contrast to the continuum of sliding canards in \esref{pinch1V0}). The two different cases are shown in \fref{fig:fil1}. 

The weak and strong canards in the unpinched system lay at $\cc{w,z}=\cc{-1/(2\eps)^\eps,\mu y/2}$ and $\cc{w,z}=\cc{-(\mu/2\eps)^\eps,y/2}$ respectively. The weak canard clearly avoids the pinch region $|w|<1$ for $\eps<1/2$, and then lies at $\cc{W,z}=\cc{1-1/(2\eps)^\eps,\mu y/2}$. The strong canard also avoids the pinch region if $\mu>2\eps$, and is given by $\cc{W,z}=\cc{1-1/(2\eps)^\eps, y/2}$. If $\mu<2\eps$ the strong canard falls inside the pinch region and is not part of \esref{pinch2} for $W\neq0$; if it exists it is part of the sliding dynamics on $W=0$. Indeed we see that is exactly the case in \fref{fig:fil1}.

\subsection{Linearizing about the weak canard}\label{sec:weak}

We complete this study by showing that the pinched approximation of the folded node possess $\langle\frac{1-\mu}{2\mu}\rangle$ secondary canards, where $\langle n\rangle$ denotes the largest integer less than $n$. These canards rotate around the weak canard near $y=0$, with rotation numbers taking all integers from $1$ to {$\langle\frac{1-\mu}{2\mu}\rangle$}, then connect to sliding solutions in the attracting and repelling sliding regions that take them to $y\rightarrow\pm\infty$, as sketched in \fref{fig:fil2ndary}. Note that the rotation takes place in $W<0$, therefore secondary canards satisfy $W\le0$. 

\begin{figure}[h!]\centering\includegraphics[width=0.5\textwidth]{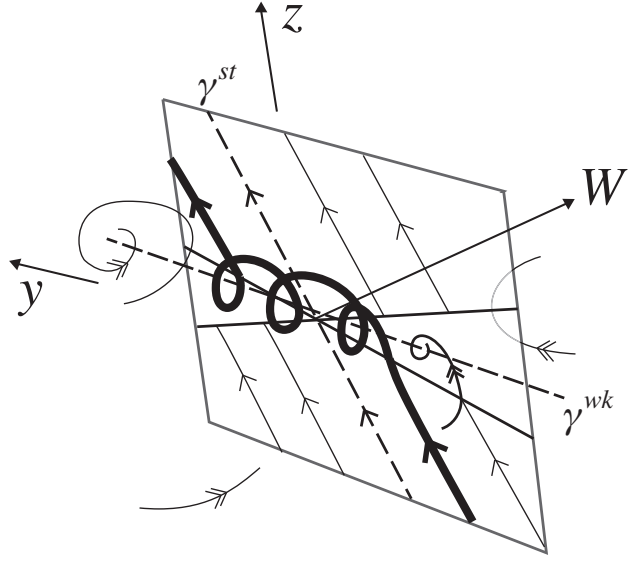}
\vspace{-0.3cm}\caption{\footnotesize\sf Primary canards (dashed) and a secondary canard (bold) with rotation number 3 in the pinched $\bb{W,y,z}$ system. }\label{fig:fil2ndary}\end{figure}

To solve \esref{pinch2} we treat the $W>0$ and $W<0$ systems separately, making different approximations in the two regions about the dominant singularities. In the region $W<0$, we linearize about the weak canard at $\cc{W,z}=\cc{1-(2\eps)^{-\eps},\mu y/2}$. In the region $W>0$, we expand about the tangency to the switching manifold at $(W,z)=\bb{0,-{\mu y}/{4\eps}}$. To leading order these give
\begin{eqnarray}\label{W<0}
{\rm in}& W>0:&\left\{\begin{array}{ccl}
\dot W&=&\frac{\mu y+4z\eps}{2\eps}+\ord{W},\\
\dot y&=&1,\\
\dot z&=&\eps+\frac{1+\mu}2+\ord{W},
\end{array}\right.\label{W>0}\\
{\rm in}& W<0:&\left\{\begin{array}{ccl}
\dot W&=&\frac{\mu y-2z}{(2\eps)^\eps}+\frac{\mu y}\eps(W-1+\frac1{(2\eps)^\eps})+\ord{\delta W^2,\eps\delta z\delta W},\\
\dot y&=&1,\\
\dot z&=&\frac\mu2+(2\eps)^{\eps-1}(W-1+\frac1{(2\eps)^\eps})+\ord{\delta W^2},
\end{array}\right.\label{W<0}
\end{eqnarray}
where $\delta W=W-1+(1/2\eps)^\eps$ and $\delta z=z-\mu y/2$. In the remainder of this section we omit the error terms, and find solutions to the truncated local equations. Note that
$1/(2\eps)^\eps\approx1-\eps\log(2\eps)+\ord{(\eps\log2\eps)^2}$
deviates quickly from unity as $\eps$ increases from zero, so we cannot approximate it by unity. Approximating around the weak canard in $W<0$ leads to a slight shift in the sliding region (\ref{Wslidingregion}). The boundary where the $W>0$ system is tangent to $W=0$ is given, as before, by $z/y=-\mu/4\eps$. The boundary where the $W<0$ system is tangent to $W=0$ is now given by 
\begin{equation}\label{Wasl}
z/y=\mu(\eps-1+(2\eps)^\eps)/2\eps\;.
\end{equation}

\begin{figure}[h!]\centering\includegraphics[width=0.8\textwidth]{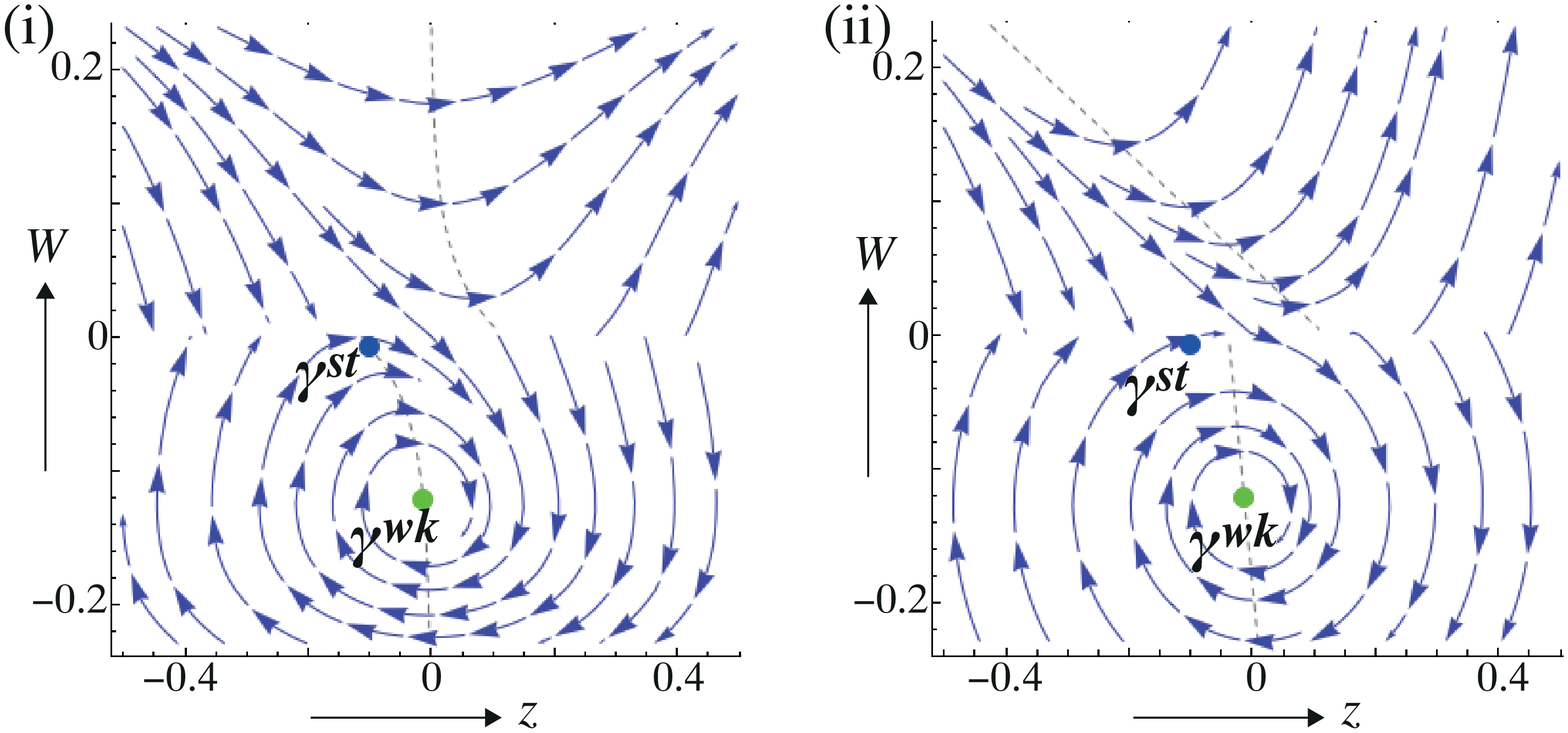}
\vspace{-0.3cm}\caption{\footnotesize\sf  The pinched system (i) and the expansion (\ref{W>0})-(\ref{W<0}). Flow simulated in the plane $y=-1$ for $\eps=0.05$, $\mu=1/8.5$. The dotted curves show the nullcline $\dot W=0$. }\label{fig:hermitefield}\end{figure}

Let us now find the canards in this approximation. We will consider only those that satisfy the following properties:
\begin{enumerate}
\item[${\cal C}1$.] $y(0)=z(0)=0$. Because canards must satisfy $W\rightarrow0$ as $t\rightarrow\infty\pm$, and the \esref{W<0}-(\ref{W>0}) is symmetric under the substitution $\cc{y,z,t}\mapsto\cc{-y,-z,-t}$, canards are expected to inherit this symmetry, implying that $y=z=0$ at $t=0$.   
\item[${\cal C}2$.] $\ddot W(t_c)<\dot W(t_c)=W(t_c)=0$ at some $t=t_c\neq0$. This is because the conditions $W=\dot W=0$ define the boundary of the sliding region (\eref{Wasl}), at which solutions can pass between $W<0$ and the invariant sliding region on $W=0$. They can only do so if the solution is curving into ($\ddot W<0$) the region $W<0$. 
\end{enumerate} 
We must then find solutions of the $W<0$ system that satisfy these two conditions. Considering \esref{W<0} we have, noting $\dot z=dz/dy=dz/dt$,  
$$\ddot z=\frac1{(2\eps)^{1-\eps}}\dot W=\frac{\mu y-2z}{2\eps}+\frac{\mu y}\eps(\dot z-\frac\mu2)\;.$$
Letting 
$z=\zeta+\mu y/2$ and $y=\tau\sqrt{\eps/\mu}$, this rearranges to the Hermite equation
\begin{equation}
\zeta''-\tau\zeta'+\frac1{\mu}\zeta=0\;,
\end{equation}
whose general solution can be written
$$\zeta(\tau)=\tau\zeta'(0)\;_1F_1\bb{\frac{\mu-1}{2\mu},\frac32,\frac{\tau^2}2}+\zeta(0)\;_1F_1\bb{-\frac1{2\mu},\frac12,\frac{\tau^2}2}\;,$$
with derivative
$$\zeta'(\tau)=\zeta'(0)\;_1F_1\bb{\frac{\mu-1}{2\mu},\frac12,\frac{\tau^2}2}-\frac{\tau\zeta(0)}\mu\;_1F_1\bb{1-\frac1{2\mu},\frac32,\frac{\tau^2}2}\;,$$
 in terms of the confluent hypergeometric function $_1F_1$ (also known as Kummer's function $M$ where $ _1F_1(\alpha,\beta;\gamma)=M(\alpha,\beta,\gamma)$, see \cite{as}). 

Applying condition ${\cal C}1$ above, at a point where $y(0)=z(0)=0$ we have $\dot W(0)=0$, hence $\ddot z(0)=0$. In the transformed coordinates this gives initial conditions
$$\zeta(0)=0,\qquad\zeta''(0)=0\;,$$
the former of which simplifies the solution above to
\begin{equation}
\frac{\zeta(\tau)}{\zeta'(0)}=\tau\;_1F_1\bb{\frac{\mu-1}{2\mu},\frac32,\frac{\tau^2}2}\;,\qquad
\frac{\zeta'(\tau)}{\zeta'(0)}=\;_1F_1\bb{\frac{\mu-1}{2\mu},\frac12,\frac{\tau^2}2}\;,
\end{equation}
or in terms of the Gamma function $\Gamma$ and Hermite polynomials $H_n$ \cite{as},
\begin{eqnarray*}
\zeta(\tau)=\mbox{\small$  \frac{i\zeta'(0)e^{\tau^2/4}\Gamma\bb{\frac{\mu+1}{2\mu}}}{\sqrt\pi\bb{i\frac{\Gamma\bb{\frac{\mu+1}{2\mu}}}{\Gamma\bb{-\frac1{2\mu}}}+\frac{\Gamma\bb{\frac{2\mu+1}{2\mu}}}{\Gamma\bb{\frac{\mu-1}{2\mu}}}}} \left\{
\frac{\Gamma\bb{\frac{2\mu+1}{2\mu}}}{\Gamma\bb{\frac{\mu-1}{2\mu}}}D^+\bb{1,\mu,i\tau} - D^+\bb{0,\mu,\tau}   \right\}   $}   \;,
\end{eqnarray*}
where
\begin{eqnarray*}
D^\pm\sq{m,\mu,\tau}=2^{\pm(\mu+2)/2\mu}e^{-\tau^2/4}H_{m\mp1/\mu}\bb{\frac{\tau}{\sqrt2}}\;.
\end{eqnarray*}
Substituting back in $\tau=t\sqrt{\mu/\eps}$, we find that in $W<0$ there exist solutions given by
\eqblock{W<0sol}{
W(t)&=&1-(2\eps)^{-\eps}(1+\eps\mu-2\eps\dot z(0)_1F_1\bb{\frac{\mu-1}{2\mu},\frac12,\frac{\mu t^2}{2\eps}})\;,\\
y(t)&=&t\;,\\
z(t)&=&\frac{\mu t}2+t\dot z(0)\;_1F_1\sq{\frac{\mu-1}{2\mu},\frac32,\frac{\mu t^2}{2\eps}}\;,
}
and these form the portions of any secondary canards that lie in $W<0$, outside the switching manifold. 
To this we must apply the second condition, ${\cal C}2$, to pick out solutions in \esref{W<0sol} that tangentially touch (or {\it graze}) the boundaries of the sliding regions on $W=0$, where they connect to sliding solutions that form the tails of the canards. Substituting the conditions $W(t_c)=\dot W(t_c)=0$ at some $t=t_c\neq0$ into \eref{W<0}, we find that $z(t_c)$ and $\dot z(t_c)$ are given by
\begin{equation}\label{ztc}
z(t_c)/\mu t_c=\dot z(t_c)+\eps(1-\mu)/2\eps=(1+\eps-(2\eps)^\eps)/{2\eps}\;.
\end{equation}
\Fref{fig:hermite} shows a simulation of solutions given by \esref{W<0sol} subject to the boundary conditions (\ref{ztc}). 

The number and geometry of the secondary canards, specifically the number of rotations they make around the weak canard, are easily found as follows. The confluent hypergeometric function $_1F_1(-a,b,c)$ has $2\langle a+1\rangle$ real zeros \cite{as} (with $\langle n\rangle$ denoting the largest integer smaller than $n$), between which the function  oscillates through $2\langle a+1\rangle-1$ maxima/minima, and between these the functions makes $\langle a\rangle$ complete oscillations. These oscillations form the rotations of the secondary canards. 

Using the boundary conditions $W(t_c)=\dot W(t_c)=0$, a given value of $t_c>0$ picks out one of the solutions (\ref{W<0sol}), with $\langle \frac{1-\mu}{2\mu}\rangle$ rotations. The boundary conditions may be satisfied at any one of the maxima, so any of the rotation numbers from $1$ up to $\langle \frac{1-\mu}{2\mu}\rangle$ are obtained, by different solutions with unique values of $t_c$. Hence there exist $\langle \frac{1-\mu}{2\mu}\rangle$ secondary canards with rotation numbers $1,2,3,...,\langle \frac{1-\mu}{2\mu}\rangle$.

\begin{figure}[h!]\centering\includegraphics[width=0.8\textwidth]{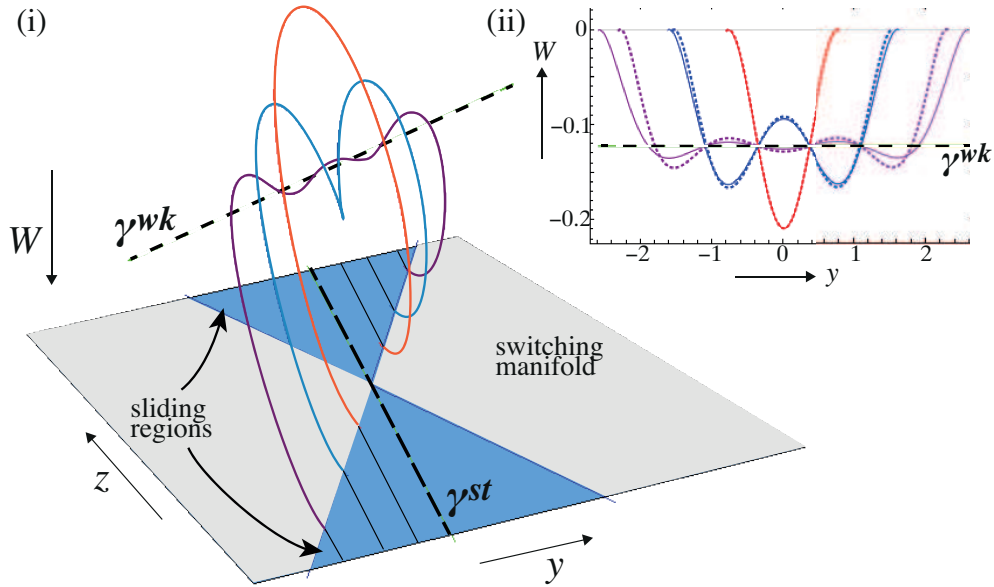}
\vspace{-0.3cm}\caption{\footnotesize\sf Secondary canards for $\eps=0.05$ and $\mu=1/8.5$, from the Hermite solution (\ref{W<0sol}). Inset: the Hermite solution (full curves) is compared to the approximation using (\ref{Happ}) (dotted curves). Long dashes indicate segments of sliding. The weak and strong canards $\gamma^{wk,st}$ are shown.}\label{fig:hermite}\end{figure}

As a final remark, the complicated expressions above can be simplified with the large parameter asymptotic approximation
\begin{equation}\label{Happ}
_1F_1\bb{a_,b_,x}=\mbox{\small$\frac{\Gamma(b)}{\sqrt\pi}\bb{\frac x2(b-2a)}^{(1-2b)/4}e^{x/2}\cos\bb{\sqrt{2x(b-2a)}+\frac\pi4(1-2b)}$}\;
\end{equation}
(see e.g.  Eq.13.5.14 of \cite{as})
This is sufficient to approximate solutions for small $x=\tau^2/2$, but this exponential approximation (compared to the exact solution in \fref{fig:hermite}(ii)) is not accurate enough to correctly give the correct number of, or number of rotations of, secondary canards.

\section{Closing remarks}\label{sec:conc}

Pinching captures the key geometry  -- singularities and bifurcations -- necessary to provide a discontinuous model of a {singularly} perturbed dynamical system. As a method of approximation it is purely qualitative, yet it appears to accurately describe singular features such as both primary and secondary canards in the folded node studied here. In a previous study of the van der Pol oscillator \cite{jd10}, pinching was also shown to capture the maximal canard, and to distinguish between a canard explosion and a Hopf bifurcation. 

Our main aim with this work is to help illuminate the bridge between smooth and piecewise-smooth models of dynamical systems, by showing how closely the phenomena of singular perturbation and discontinuity are related. Canards in slow-fast systems have now been known for more than three decades, having first been {studied using non-standard analysis \cite{b81}, and later with standard {tools such as matched asymptotic expansions and} geometric singular {perturbation} theory \cite{e83,j95}. Canards in discontinuous system lay unannounced in the work of Filippov \cite{f64,f88} for considerably longer. Only recently has the link between folds in critical manifolds of slow-fast systems, and the two-fold singularity in discontinuous systems, become clear, through the methods of regularization \cite{tls09} and of pinching \cite{jd10}, the former applying a topological equivalence in the singular limit $\eps=0$, the latter approximating the geometry for $\eps\neq0$.


\bigskip\noindent{\bf Acknowledgements}. MD acknowledges the support of EPSRC through grant EP/E032249/1 and the Department of Engineering Mathematics at the University of Bristol (UK) where part of this work was completed; he also acknowledge the support of the INRIA large-scale initiative \href{https://www.rocq.inria.fr/sisyphe/reglo/regate.html}{REGATE (Regulation of the GonAdoTropE axis)}. MRJ's research is supported by EPSRC Grant Ref: EP/J001317/1.

\bibliography{../grazcat}

\begin{thebibliography}{10}

\bibitem{as}
M.~Abramowitz and I.~Stegun.
\newblock {\em Handbook of Mathematical Functions}.
\newblock Dover, 1964.

\bibitem{russian2}
M.~A. Aizerman and F.~R. Gantmakher.
\newblock On the stability of equilibrium positions in discontinuous systems.
\newblock {\em Journal of Applied Mathematics and Mechanics (translated from
  the Russian Ob ustoichivosti polozhenii ravnovesiia v razryvnykh sistemakh)},
  24:283--93, 1960.

\bibitem{ben90}
E.~Beno\^{i}t.
\newblock Canards et enlacements.
\newblock {\em Publications math\'ematiques de l'IH\'ES}, 72:63--91, 1990.

\bibitem{b81}
E.~Beno\^{i}t, J.~L. Callot, F.~Diener, and M.~Diener.
\newblock Chasse au canard.
\newblock {\em Collect. Math.}, 31-32:37--119, 1981.

\bibitem{berry94red}
M.~V. Berry.
\newblock Asymptotics, singularities and the reduction of theories.
\newblock {\em Logic, Methodology and Philosophy of Science IX}, pages
  597--607, 1994.

\bibitem{brons06}
M.~Brons, M.~Krupa, and M.~Wechselberger.
\newblock Mixed mode oscillations due to the generalized canard phenomenon.
\newblock {\em Fields Institute Communications}, 49:39--63, 2006.

\bibitem{cj09}
A.~Colombo, M.~di~Bernardo, E.~Fossas, and M.~R. Jeffrey.
\newblock Teixeira singularities in 3{D} switched feedback control systems.
\newblock {\em Systems and Control Letters}, 59(10):615--622, 2010.

\bibitem{des12}
M.~Desroches, J.~Guckenheimer, B.~Krauskopf, C.~Kuehn, H.~M. Osinga, and
  M.~Wechselberger.
\newblock Mixed-mode oscillations with multiple time scales.
\newblock {\em SIAM Rev.}, 54(2):211--288, 2012.

\bibitem{jd10}
M.~Desroches and M.~R. Jeffrey.
\newblock Canards and curvature: nonsmooth approximation by pinching.
\newblock {\em Nonlinearity}, 24:1655--1682, 2011.

\bibitem{d08}
M.~Desroches, B.~Krauskopf, and H.~M. Osinga.
\newblock Mixed-mode oscillations and slow manifolds in the self-coupled
  {F}itzhugh-{N}agumo system.
\newblock {\em Chaos}, 18(1):015107, 2008.

\bibitem{des10}
M.~Desroches, B.~Krauskopf, and H.~M. Osinga.
\newblock Numerical continuation of canard orbits in slow-fast dynamical
  systems.
\newblock {\em Nonlinearity}, 23(3):739--765, 2010.

\bibitem{bc08}
M.~di~Bernardo, C.~J. Budd, A.~R. Champneys, and P.~Kowalczyk.
\newblock {\em Piecewise-Smooth Dynamical Systems: Theory and Applications}.
\newblock Springer, 2008.

\bibitem{bkn02}
M.~di~Bernardo, P.~Kowalczyk, and A.~Nordmark.
\newblock Bifurcations of dynamical systems with sliding: derivation of
  normal-form mappings.
\newblock {\em Physica D}, 170:175--205, 2002.

\bibitem{e83}
W~Eckhaus.
\newblock Relaxation oscillations including a standard chase on {F}rench ducks.
\newblock {\em Lect. Notes Math.}, 985:449--494, 1983.

\bibitem{f79}
N.~Fenichel.
\newblock Geometric singular perturbation theory.
\newblock {\em J. Differ. Equ.}, 31:53--98, 1979.

\bibitem{f64}
A.~F. Filippov.
\newblock Differential equations with discontinuous right-hand side.
\newblock {\em American Mathematical Society Translations, Series 2},
  42:19--231, 1964.

\bibitem{f88}
A.~F. Filippov.
\newblock {\em Differential Equations with Discontinuous Righthand Sides}.
\newblock Kluwer Academic Publ. Dortrecht, 1988.

\bibitem{guck00}
J.~Guckenheimer, K.~Hoffman, and W.~Weckesser.
\newblock Numerical computation of canards.
\newblock {\em Int. J. Bifurc. Chaos}, 10(12):2269--2687, 2000.

\bibitem{j11prl}
M.~R. Jeffrey.
\newblock Non-determinism in the limit of nonsmooth dynamics.
\newblock {\em Physical Review Letters}, 106(254103):1--4, 2011.

\bibitem{jh09}
M.~R. Jeffrey and S.~J. Hogan.
\newblock The geometry of generic sliding bifurcations.
\newblock {\em SIAM Review}, 53(3):505--525, 2011.

\bibitem{j95}
C.~K. R.~T. Jones.
\newblock {\em Geometric singular perturbation theory}, volume 1609 of {\em
  Lecture Notes in Math. pp. 44-120}.
\newblock Springer-Verlag (New York), 1995.

\bibitem{krupa08}
M.~Krupa, N.~Popovi\'{c}, and N.~Kopell.
\newblock Mixed-mode oscillations in three time-scale systems: A prototypical
  example.
\newblock {\em SIAM J. Appl. Dyn. Syst.}, 7(2):361--402, 2008.

\bibitem{tls09}
J.~Llibre, P.~R. da~Silva, and M.~A. Teixeira.
\newblock Study of singularities in nonsmooth dynamical systems via singular
  perturbation.
\newblock {\em SIAM J. App. Dyn. Sys.}, 8(1):508--526, 2009.

\bibitem{russian1}
Yu.~I. Neimark and S.~D. Kinyapin.
\newblock On the equilibrium state on a surface of discontinuity.
\newblock {\em Radiophysics and Quantum Electronics (translated from the
  Russian Izvestiia vysshikh uchebnykh zavedenii. Radiofizika.)}, 3:694--705,
  1960.

\bibitem{pop08}
N.~Popovi\'{c}.
\newblock Mixed-mode dynamics and the canard phenomenon: Towards a
  classification.
\newblock {\em J. Phys.: Conf. Ser.}, 138(012020), 2008.

\bibitem{sw01}
P.~Szmolyan and M.~Wechselberger.
\newblock Canards in $\mathbb{R}^3$.
\newblock {\em J. Differ. Equ.}, 177:419--453, 2001.

\bibitem{t93}
M.~A. Teixeira.
\newblock Generic bifurcation of sliding vector fields.
\newblock {\em J.Math.Anal.Appl.}, 176:436--457, 1993.

\bibitem{w05}
M.~Wechselberger.
\newblock Existence and bifurcation of canards in $\mathbb{R}^3$ in the case of
  a folded node.
\newblock {\em SIAM J. App. Dyn. Sys.}, 4(1):101--139, 2005.

\bibitem{w12}
M.~Wechselberger.
\newblock A propos de canards (apropos canards).
\newblock {\em Trans. Amer. Math. Soc}, 364:3289--3309, 2012.

\end{thebibliography}
\bibliographystyle{plain}

\end{document}